\documentclass[12pt]{article}
\usepackage{amsfonts}
\begin{document}
\title{\bf   Orbital stability of standing waves of Two-component
Bose-Einstein condensates with internal atomic Josephson Junction
\author{ H. Hajaiej\\
{\footnotesize{\em  College of Sciences, King Saud University, 11451 Riyadh}}\\
{\footnotesize{\em Department of Mathematics}}\\
 {\footnotesize{\em E.mail : hhajaiej@ksu.edu.sa}} }}
\date{}
\maketitle
\begin{abstract}
In this paper, we prove existence, symmetry and uniqueness of
standing waves for a coupled Gross-Pitaevskii equations modeling
component Bose-Einstein condensates(BEC) with an internal atomic
Josephson junction. We will then address the orbital stability of
these standing waves and characterize their orbit.
\end{abstract}
\section{Introduction}

The dynamics of a model of a two-component (BEC) irradiated by an
external electromagnetic field are given by the following two
coupled nonlinear Schr\"odinger equations :
$$\left.\begin{array}{ll}
i\partial_t \psi_j &= - \displaystyle{\frac{1}{2}} \Delta \psi_j +
 \displaystyle{\frac{\gamma^2}{2}}|x|^2\psi_j + \beta_{jj}|\psi_j|^2\psi_j
 + \beta_{ji}|\psi_i|^2 \psi_j+\lambda\psi_i+\delta\psi_j\\
 \psi_j(0,x) &= \psi^0_j\end{array}\right\}\eqno{(1.1)}$$
 $i\neq j = (1,2), (t,x) \in \mathbb{R} \times \mathbb{R}^N, N = 1, 2,
 3$.\\
 $V(x) = \displaystyle{\frac{\gamma^2}{2}}|x|^2$ is the trapping
 potential, $\gamma > 0$.\\
 $\beta_{12} = B_{21}$ is the inter-specific scattering length,
 while $\beta_{11}$ and $\beta_{22}$ are the intra ones. $\lambda$
 is the rabi frequency related to the external electric field. It
 is the effective frequency to realize the internal atomic Josephson
 junction by a Raman transition, $\delta$ is the detuning constant
 for the Raman transition.
 (1.1) arises in modelling  BEC composed of atoms in two hyperfine
 states in the same harmonic map [1]. Recently, BEC with multiple
 species have been realized in experiments, ([2] and references
 therein) and many interesting phenomena, which do not appear in
 the single component BEC, have been observed  in the
 multi-component BEC.
  The simplest multi-component BEC can be viewed
 as a binary mixture, which can be used as  a model to produce
 atomic lazer. To our knowledge, the first experiment in this
 framework was done quite recently [9], this has opened the way to
 many other groups of research who carried out the study of such
 problems for two-component BEC theoretically and experimentally.\\
 In this paper, we consider a binary BEC model in which there is an
 irradiation with an electromagnetic field, this causes a
 Josephson-type oscillation between the two species. These
 condensates are extremely important in physics and nonlinear optics
 since it is possible to measure the relative phase of one component
 with respect to the other one [Lemma 2.1, 2]. Controlling the
 relative phase, it is also possible to produce vortices, [7], [9].
 A more detailed account is given in [1].

 A standing wave for (1.1) is a function $(\psi_1, \psi_2) = (e^{-i\mu_1t}\Phi_1 ,
 e^{-i\mu_2t}\Phi_2)$ solving this NLS. Thus it satisfies the
 following $2 \times 2(\mathbb{C})$ elliptic system :
 $$\left\{ \begin{array}{ll}
 \mu_1 \Phi_1 &= \Big[- \displaystyle{\frac{1}{2} \Delta +\frac{\gamma^2}{2}}
 |x|^2 + \delta + \beta_{11}|\Phi_1|^2+ \beta_{12}|\Phi_2|^2\Big]\Phi_1 +
 \lambda \Phi_2\\
 \mu_2 \Phi_2 &= \Big[- \displaystyle{\frac{1}{2} \Delta + \frac{\gamma^2}{2}}
 |x|^2 + \delta+ \beta_{22} |\Phi_2|^2 + \beta_{12} |\Phi_1|^2\Big]
 \Phi_2+ \lambda \Phi_1
 \end{array}\right.\eqno{(1.2)}$$
 Ground state solutions of (1.2) are the minimizes of the following
 constrained variational problem :\\
 For two prescribed real numbers $c_1$ and $c_2$
 $$\hat{I}_{c_1,c_2} = \inf_{(\Psi_1,\Psi_2) \in \hat{S}_{c_1,c_2}}
 \hat{E} (\Psi_1,\Psi_2)\eqno{(1.3)}$$
 $$\hat{S}_{c_1,c_2} = \left\{ (\Psi_1,\Psi_2) \in \Sigma_{\mathbb{C}}
 (\mathbb{R}^N) \times \Sigma_{\mathbb{C}}(\mathbb{R}^N)
 \times  : \int|\Psi_1|^2 = c^2_1 \mbox{ and }
 \int |\Psi_2|^2 = c^2_2\right\}.$$
 $$\sum (\mathbb{R}^N) = \left\{ u \in H^1(\mathbb{R}^N) : \int_{\mathbb{R}^N}
 |x|^2 u^2(x)dx < \infty\right\}$$
 $$|u|^2_{\sum(\mathbb{R}^N)} = |u|^2_2 + |\nabla u|^2_2 + ||x|u|^2_2$$
 $$\Sigma_{\mathbb{C}}(\mathbb{R}^N) =
 \left\{ z = (u,v) \simeq u + iv : (u,v) \in \sum (\mathbb{R}^N)
 \times \sum (\mathbb{R}^N)\right\}$$
 $$\|z\|^2_{\sum_{\mathbb{C}}(\mathbb{R}^N)} = \|z\|^2_2
 + \|\nabla z\|^2_2 + \||x| z\|^2_2.$$
 $$\hat{E}(\Psi) = \hat{E}_0 (\Psi_1,\Psi_2)+2\lambda\int Re(\Psi_1\overline{\Psi_2})dx,$$
 with $\bar{f}$ denoting the conjugate part of $f$ and Ref its real one.
 $$\begin{array}{ll}
 \hat{E}_0(\psi_) &= \hat{E}_0(\Psi_1,\Psi_2) =
 \displaystyle{\int_{\mathbb{R}^N}\frac{1}{2}}
 \Big[|\nabla \Psi_1|^2_2 + |\nabla \Psi_2|^2_2 +
 \gamma^2|x|^2(|\Psi_1|^2 +|\Psi_2|^2)\Big]\\
 &+ \delta|\Psi_1|^2+ \displaystyle{\frac{1}{2}}\beta_{11}
 |\psi_1|^4+
 \displaystyle{\frac{1}{2}}\beta_{22}|\Psi_2|^4+
 \beta_{12}|\Psi_1|^2|\Psi_2|^2\Big\}dx
 \end{array}\eqno{(1.5)}$$
 As proved in iii) and iv) of Lemma 2.1 of [2], solving the
 constrained  minimization problem (1.3) is equivalent to study
 auxiliary minimization problem :
 $$\widetilde{\hat{I}}_{c_1,c_2} = \inf_{(\Psi_1,\Psi_2) \in \hat{S}_{c_1,c_2}}
 \widetilde{\hat{E}}(\Psi_1,\Psi_2)\eqno{(1.6)}$$
 where
 $$\widetilde{\hat{E}}(\Psi) = \widetilde{\hat{E}}(\Psi_1,\Psi_2) =
 \hat{E}_0(\Psi)- 2|\lambda| \int_{\mathbb{R}^N}|\Psi_1||\Psi_2|dx\eqno{(1.7)}$$
 The main objective of the present work is to show the orbital
 stability of standing waves of (1.1). To reach this goal, we
 will first solve (1.6) for real-valued functions :
 $$\widetilde{I}_{c_1,c_2) \in S_{c_1,c_2}}
 \widetilde{E}(u_1,u_2) = \widetilde{E}(u).\eqno{(1.8)}$$
 $$\widetilde{E}(u)= \widetilde{E}(u_1,u_2) = \widetilde{E}_0(u_1,u_2)
 - 2|\lambda| \int|u_1||u_2|dx.\eqno{(1.9)}$$
 $$S_{c_1,c_2} = \left\{(u_1,u_2) \in \sum(\mathbb{R}^N)
 \times \sum(\mathbb{R}^N) : \int_{\mathbb{R}^N}
 u^2_1(x) = c^2_1 \mbox{ and } \int_{\mathbb{R}^N}
 u^2_2(x) = c^2_2\right\}.$$
 We will first prove existence, symmetry uniqueness of minimizers
 of (1.8). Then we will use these qualitative properties to solve
 the constrained variational problem (1.6), which is in itself a key
 step to show the orbital stability of standing waves  and to
 characterize their  orbit.\\
 In this paper, we will concentrate our study on the critical case
 $N = 2$, which is,  from the mathematical point of view,   the most
 challenging case.\\
 Our paper is organized as follows. In section 2, we will give some
 important definitions and preliminary results. Then  we will derive
 some qualitative  properties of the energy functional and the
 minimization problem, this will be the key ingredient to study the
 orbital stability of standing waves in the last section.\\
 We will focus our study on the case $N = 2$ but we will give clear
 and complete indications about $N=1$ and $N =3$.
 \section{Notation,  Definitions and Preliminary Results}
 $H^1(\mathbb{R}^N)$ is the usual Hilbert space
 $$\sum(\mathbb{R}^N) = \{u \in H^1(\mathbb{R}^N) :
 \int_{\mathbb{R}^N}|x|^2|u|^2 dx < \infty\}$$
 $$|u|^2_{\sum(\mathbb{R}^N)} = |u|^2_2 + |\nabla u|^2_2 + ||xu|^2_2$$
 $|u|_p$ is the standard norm of the $L^p(\mathbb{R}^N)$ space
 $$H^1(\mathbb{R}^N,\mathbb{C}) = \{z = (u,v) \in H^1(\mathbb{R}^N) \times
 H^1(\mathbb{R}^N)\}.$$
 We shall identify $z = (u,v)$ with $u+iv \in
 H^1(\mathbb{R}^N,\mathbb{C})$.\\
 For $z \in H^1(\mathbb{R}^N,\mathbb{C}) $, $\|z\|^2_{H^1(\mathbb{R}^N,\mathbb{C})}
 = \|z\|^2_2 + \|\nabla z\|^2_2$
 $$\|z\|^2_2 = |u|^2_2 + |v|^2_2 \mbox{ and } \|\nabla z\|^2_2
  = |\nabla u|^2_2+ |\nabla v|^2_2.$$
 Here and elsewhere $|\;|_q$ denotes the usual norm in
 $L^q(\mathbb{R}^N)$ and $\|\;\|_q$ is the standard norm in
 $L^q(\mathbb{R}^N,\mathbb{C})$.
 $$\Sigma_{\mathbb{C}}(\mathbb{R}^N) = \left\{z \in H^1(\mathbb{R}^N,\mathbb{C})
 \int |x|^2|z|^2 dx < \infty\right\}.$$
 $\sum(\mathbb{R}^N) \times \sum (\mathbb{R}^N)$ and $\sum_{\mathbb{C}}(\mathbb{R}^N)
 \times \sum_{\mathbb{C}}(\mathbb{R}^N)$ are equipped with the
 standard cartesian norms.\\
 For fixed real numbers $c_1$ and $c_2$, we define
 $$Z_{c_1,c_2} = \{(z_1,z_2) \in \hat{S}_{c_1,c_2} :
 \widetilde{\hat{E}}(z_1,z_2) = \widetilde{\hat{I}}_{c_1,c_2}\}\eqno{(2.1)}
$$
$$W_{c_1,c_2} = \{(u_1,u_2) \in S_{c_1,c_2}, u_1 \mbox{ and }
u_2 > 0 \mbox{ and } \widetilde{E}(u_1,u_2) =
\widetilde{I}_{c_1,c_2}\}.$$ We say that $Z_{c_1,c_2}$ is stable if
:
$$Z_{c_1,c_2} \neq \emptyset \quad \mbox{ and }
\;\forall\; w = (w_1,w_2) \in Z_{c_1,c_2},$$ $\forall\; \varepsilon
> 0, \;\exists\; \delta > 0$ such that for any $\psi_0 = (\Psi^1_0,
\Psi^2_0) \in \sum_{\mathbb{C}} (\mathbb{R}^N) \times
\sum_{\mathbb{C}}(\mathbb{R}^N)$
$$\left\{ \begin{array}{ll}
\mbox{ satisfying } \|\Psi_0-w\|_{\sum_{\mathbb{C}(\mathbb{R}^N)}} <
\delta\\
\displaystyle{\inf_{z \in Z_{c_1,c_2}}}\|\Psi(t,.)
-z\|_{\Sigma_{\mathbb{R}}\times \Sigma_{\mathbb{C}}(\mathbb{R}^n)} <
\varepsilon
\end{array}\right.\eqno{(2.2)}$$
for all $t \in \mathbb{R}$, where $\Psi(t,.)$ is the unique solution
of (1.1) corresponding to the initial condition $\Psi_0$.  (Note
that in [1], the authors have solved the Cauchy problem (1.1)
under the assumptions of Lemma 2.3 below).\\
{\bf Lemma 2.1.} [Lemma 4.4,8]

$\sum(\mathbb{R}^N)$ is compactly embedded  in $L^q(\mathbb{R}^N)$
for any $q$ such that $q < \frac{2N}{N-2}$.\\
{\bf Lemma 2.2} Let $N = 2$,
$$(A_1) \left\{ \begin{array}{l}
\beta_{ij} < 0, 1 \leq i,j \leq 2 \mbox{ and }\\
\beta_{11} c^2_1 + \beta_{12} c_1c_2 > - c_b\\
\beta_{22} c^2_2 + \beta_{12} c_1c_2 > - c_b
\end{array}\right.$$
$c_b$ is defined as the best constant in the Gagliardo-Nirenberg
inequality
$$\int_{\mathbb{R}^2}u^4 \leq \frac{1}{c_b} |\nabla u|^2_2 |u|^2_2\eqno{(2.3)}$$
Then :
\begin{enumerate}
\item The minimization problem (1.8) is well-posed and any
minimizing sequence of (1.8) is bounded in $\sum(\mathbb{R}^2)
\times \sum (\mathbb{R}^2)$.
\item Any minimizing sequence of (1.8) is relatively compact in
$\sum (\mathbb{R}^2) \times \sum (\mathbb{R}^2)$, i.e, $\forall\;
u_n = (u_{n,1}, u_{n,2}) \subset S_{c_1,c_2}$ such that
$E(u_{n,1},u_{n,2}) \rightarrow I_{c_1,c_2}$, then there exists $u =
(u_1,u_2) \in \sum (\mathbb{R}^2) \times \sum(\mathbb{R}^2)$ such
that $u_n \rightarrow u$ in $\sum (\mathbb{R}^2) \times \sum
(\mathbb{R}^2)$ (up to a subsequence).
\item The functionals $\widetilde{E}$ and $\widetilde{\hat{E}}$ are
$C^1$ in $\sum(\mathbb{R}^2) \times \sum (\mathbb{R}^2)$ (resp.
$\sum_{\mathbb{C}}(\mathbb{R}^2) \times
\sum_{\mathbb{C}}(\mathbb{R}^2)$
\item $(c_1,c_2) \rightarrow \widetilde{I}_{c_1,c_2}$ is continuous
.
\end{enumerate}
{\bf Proof}
\begin{enumerate}
\item Let $(u_1,u_2) \in S_{c_1,c_2}$.\\
First using Gagliardo-Nirenberg inequality, we know that
$$\int_{\mathbb{R}^2}
|u_1|^4 \leq \frac{1}{c_b} |\nabla u_1|^2_2 |u_1|^2_2 =
\frac{1}{c_b} |\nabla u_1|^2_2 c^2_1\eqno{(2.4)}$$ and
$$\int_{\mathbb{R}^2}
|u_2|^4 \leq \frac{1}{c_b} |\nabla u_2|^2_2 |u_2|^2_2 =
\frac{1}{c_b}
  |\nabla u_2|^2_2 c^2_2.$$
  On the other hand, by Hardy inequality, we have that :
  $$\int_{\mathbb{R}^2}|u_1|^2 |u_2|^2 \leq \left(\int_{\mathbb{R}^2}
  |u_1|^4\right)^{1/2} \left(\int_{\mathbb{R}^2}|u_2|^4\right)^{1/2}
  \leq \frac{c_1c_2}{c_b} |\nabla u_1|_2|\nabla u_2|_2\eqno{(2.5)}$$
\end{enumerate}
It follows by Young inequality that :
$$\int_{\mathbb{R}^2} |u_1|^2|u_2|^2 \leq \frac{1}{2c_b}
c_1c_2 [|\nabla u_1|^2_2 + |\nabla u_2|^2_2]\eqno{(2.6)}$$ On the
other hand, we can easily prove that
$$-2 |\lambda| \int_{\mathbb{R}^2} |u_1||u_2| \geq -2|\lambda| c_1c_2\eqno{(2.7)}$$
Combining (2.4) to (2.7), we get :
$$\begin{array}{ll}
\widetilde{E}(u) = \widetilde{E}(u_1,u_2) &\geq |\nabla u_1|^2_2
\{\displaystyle{\frac{1}{2} + \frac{1}{2} \beta_{11}
\frac{c^2_1}{c_b} + \frac{1}{2} \beta_{12} \frac{c_1c_2}{c_b}}\}\\
&+ |\nabla u_2|^2_2 \{1/2 + \frac{1}{2} \beta_{22} \frac{c^2_2}{c_b}
+ \frac{1}{2} \beta_{12} \frac{c_1c_2}{c_b}\} -
|\delta||u_1|^2_2\\
&- 2|\lambda|c_1c_2 +
\displaystyle{\frac{\gamma^2}{2}\int_{\mathbb{R}^2}}|x|^2(|u_1|^2 +
|u_2|^2).
\end{array}\eqno{(2.8)}$$
$(A_1)$ enables us to conclude that the energy functional
$\widetilde{E}$ is bounded from below in $\sum(\mathbb{R}^N) \times
\sum (\mathbb{R}^N)$.\\
{\bf Remark 1}
\begin{itemize}
\item[a)] If there exists $\beta_{ij} \geq 0$, then 1) still  holds
true by replacing $\beta_j \geq 0$ by 0 in the assumption $(A_1)$.
\item[b)] In [2], the boundedness from below of the energy
functional $\widetilde{E}$ has been proved differently (page 56,
line 9).\\
More precisely : Combining Cauchy and Gagliardo Nirenberg
inequalities, the authors have  proved that
$$\begin{array}{ll}
\displaystyle{\int_{\mathbb{R}^2}} \beta_{11} |u_1|^4 &+
\beta_{22}|u_2|^4+2\beta_{12} |u_1|^2 |u_2|^2dx\geq\\
&- c_b \displaystyle{\int_{\mathbb{R}^2}}(\sqrt{|u_1|^2+|u_2|^2})^4
dx \geq -
\displaystyle{\int_{\mathbb{R}^2}}(\sqrt{|u_1|^2+|u_2|^2})^2dx
\displaystyle{\int_{\mathbb{R}^2}}
(\nabla \sqrt{|u_1|^2+|u_2|^2})^2\\
& \geq - \displaystyle{\int_{\mathbb{R}^2}}|\nabla u_1|^2 + |\nabla
u_2|^2.\hspace*{8cm} (2.8')
\end{array}$$
provided that
$$ \left.\begin{array}{ll}
\beta_{11} &> - c_b\\
\beta_{22} &> -c_b\\
\beta_{12} &\geq -c_b - \sqrt{\beta_{11}+c_b}
\sqrt{\beta_{22}+c_b}\\
\mbox{ and }\\
c^2_1 &+ c^2_2 < 1.
\end{array}\right\}.\eqno{(A_1')}$$
It seems that if one uses their approach, it is necessary to impose
the very restrictive condition : $c^2_1+c^2_2 < 1$.\\
Nevertheless the two approaches are equivalent if one considers the
same one-constrained minimization problem (1.8) with "their"
constraint $\displaystyle{\int}u^2_1+ u^2_2 = 1$.
\item[c)] The case $N = 1$ is immediate since the Gagliardo
Nirenberg inequality is not critical there.

However for $N=3$, our approach only applies when all the constants
$\beta_{ij}$ are positive. When all of then are strictly negative,
contrary to what was stated in ([2], Theorem 2.2), the minimization
problem is ill-posed. Nevertheless if $c^2_1+c^2_2 < 1$ then using
the same approach developed in [2] ((2.8)') we can easily prove that
1) and 2) still hold true if we have the following assumption :
$$\left. \begin{array}{l}
\beta_{11} > 0\\
\beta_{22} > 0\\
\beta_{11} \beta_{22} - \beta^2_{12} > 0.
\end{array}\right\}\eqno{(A_1)_{N=3}}$$
$$\beta_{11} > 0, \beta_{12} > 0 \mbox{ and } \beta_{22} > 0.\eqno{(A'_1)_{N=3}}$$
\item[d)] Let us finally emphasize that all the results of this
section hold true provided that the constrained minimization problem
is well-posed.
\end{itemize}
{\bf Proof of 2)}

By 1), we can conclude that any minimizing sequence $u_n =
(u_{n,1},u_{n,2})$ of (1.8) is bounded in $\sum (\mathbb{R}^2)
\times \sum (\mathbb{R}^2)$. Therefore up to a subsequence (that we
will also denote by $(u_n))$, there exists $(u_1,u_2) \in \sum
(\mathbb{R}^2) \times \sum(\mathbb{R}^2)$ such that $u_{n,1}
\rightharpoonup u_1$ and $u_{n,2} \rightharpoonup u_2$ in
$\sum(\mathbb{R}^2)$.\\
By Lemma 2.1, $u_{n,1} \rightarrow u_1$ and $u_{n,2} \rightarrow
u_2$ in $L^p(\mathbb{R}^2),\;\; \forall\; 2 \leq p <
\infty.\hfill{(2.9)}$\\
 Now note that by the lower semi-continuity
of the norm $|\;|_{\sum(\mathbb{R}^N)}$, we certainly have :
$$\begin{array}{ll}
\frac{1}{2} |\nabla u_1|^2_2 &+ \frac{1}{2} |\nabla u_2|^2_2 +
\frac{\gamma^2}{2} \displaystyle{\int_{\mathbb{R}^2}}|x|^2
(u_1|^2+|u_2|^2)dx\\
&\leq \liminf\Big(\frac{1}{2} |\nabla u_{n,1}|^2_2 +
\frac{1}{2}|\nabla u_{n,2}|^2_2 + \frac{\gamma^2}{2}
\displaystyle{\int_{\mathbb{R}^2}}|x|^2 +
|u_{n,1}|^2+|u_{n,2}|^2)dx\Big)
\end{array}\eqno{(2.10)}$$
On the other hand, by (2.9), we have that :
$$\int_{\mathbb{R}^2}|u_{n,1}|^4 \rightarrow \int|u_1|^4$$
$$\int_{\mathbb{R}^2}|u_{n,2}|^4 \rightarrow \int|u_2|^4\eqno{(2.11)}$$
Thus using the dominated convergence theorem, we can deduce that.
$$\int_{\mathbb{R}^2}|u_{n,1}|^2 |u_{n,2}|^2 \rightarrow \int_{\mathbb{R}^2}|u_1|^2
|u_2|^2\eqno{(2.12)}$$
 Indeed since $u_n \rightarrow u$ in
$L^4(\mathbb{R}^2) \times L^4(\mathbb{R}^2)$, there exist a
subsequence $(u_{nj,1} \subset L^4(\mathbb{R}^2)$ and a function $h
\in L^4(\mathbb{R}^2)$ such that $u_{nj,1} \rightarrow u_1$ almost
every where with $|u_{nj,1}|\leq h$.\\
Similarly, we can find $(u_{nj,2}$ and $k \subset L^4(\mathbb{R}^2)$
such that  $u_{n,2}\rightarrow u_2$ a.e with $|u_{nj,2}| \leq k)$
$$\int_{\mathbb{R}^2} |u_{nj,1}|^2 |u_{nj,2}|^2 \leq \int_{\mathbb{R}^2}
h^2k^2 dx \leq(\int h^4)^{1/2}(\int k^4)^{1/2} < \infty .$$
Therefore
$$\lim_{n\rightarrow \infty}\int_{\mathbb{R}^2}|u_{nj,1}|^2
|u_{nj,2}|^2 \leq \int_{\mathbb{R}^2} h^2k^2 dx \leq (\int
h^4)^{1/2} (\int k^4)^{1/2} < \infty$$
In the same manner, we can
prove that $\displaystyle{\lim_{n\rightarrow
\infty}\int|u_{n,1}||u_{n,2}| =
\int}|u_1||u_2|$.\hfill (2.13)\\
Combining (2.10) to (2.13), we obtain :
$$\widetilde{E}(u) = \widetilde{E} (u_1,u_2) \leq
\liminf \widetilde{E}(u_{n,1},u_{n,2}) =
\widetilde{I}_{c_1,c_2}\eqno{(2.14)}$$ But
$$\int u^2_1 = \lim_{n\rightarrow \infty} \int_{\mathbb{R}^2}u^2_{n,1} = c^2_1\quad
\mbox{ and } \int u^2_2 = \lim_{n\rightarrow
\infty}\int_{\mathbb{R}^2} u^2_{n,2} = c^2_2.$$ Thus $u = (u_1,u_2)
\in S_{c_1,c_2}$ with $\widetilde{E}(u) = \widetilde{E}(u_1,u_2) =
\widetilde{I}_{c_1,c_2}$.\\
{\bf Remark 2} :

$\bullet\quad$ In part 2) of the Lemma, we have also proved that any
minimizing sequence of (1.8) is relatively compact in
$\sum(\mathbb{R}^N) \times \sum(\mathbb{R}^N)$

$\bullet\quad$ The proofs of 3) and 4) goes exactly in the same way
as in [Proposition 3.2, 4].\\
{\bf Lemma 2.3} Under $(A_1)$, all the minimizers of (1.3) are
non-negative radial and radially decreasing.\\
{\bf Proof} First note that $\widetilde{E}(|u_1|,|u_2|) \leq
\widetilde{E}(u_1,u_2)$ for any $(u_1,u_2)$ for any $(u_1,u_2) \in
\sum (\mathbb{R}^N) \times \sum (\mathbb{R}^N)$. Therefore, we can
suppose without less of generality that $u_1$ and $u_2$ are
non-negative.\\
On the other hand, using rearrangement inequalities, [5], we know
that for any $f,g$ non-negative $\in \sum (\mathbb{R}^N)$, we have :
$$\int f^2 = \int(f^\ast)^2$$
$$\int f^4 = \int(f^\ast)^4$$
$$\int f g \leq \int f^\ast g^\ast$$
$$\int f^2 g^2 \leq \int(f^\ast)^2(g^\ast)^2$$
$$\int|x|^2 (f^\ast)^2 < \int |x|^2 f^2,$$
and
$$|\nabla f^\ast|_2 \leq |\nabla f|_2.$$
{\bf Lemma 2.4} If
$$(A_2) \left\{ \begin{array}{ll}
\beta = \left(\begin{array}{ll} \beta_{11} &\beta_{12}\\
\beta_{12} &\beta_{22}
\end{array}\right) \mbox{ is positive semi-definite }\\
\mbox{ and at lesst } :
\begin{array}{ll}
\beta_{11} - \beta_{22}
&\neq 0 \mbox{ or }\\
\beta_{11}- \beta_{12} &\neq 0 \mbox{ or }\\
\delta &\neq 0 \mbox{ or }\\
\lambda &\neq 0
\end{array}
\end{array}\right.$$
Then (1.8) has a unique minimizer.\\
{\bf Proof} [Lemma 2.2, 2]
\section{Orbital stability of standing waves of (1.1)}
In all this section, we assume that $(A_1)$ and $(A_2)$ hold true.\\
{\bf Theorem 3.1}
\begin{enumerate}
\item For any $c_1,c_2 ; \widetilde{I}_{c_1,c_2} =
\widetilde{\hat{I}}_{c_1,c_2} , Z_{c_1,c_2}\neq \emptyset$ and
$Z_{c_1,c_2}$ is stable
\item For any $z = (z_1,z_2) \in Z_{c_1,c_2}, |z| = (|z_1|,|z_2|) \in
W_{c_1,c_2}$ and $$Z_{c_1,c_2} = \left\{(e^{i\theta_1w_2},
e^{i\theta_2w_2}), (\theta_1,\theta_2) \in \mathbb{R}^2\right\},$$
where $(w_1,w_2)$ is the unique solution of (1.8).
\end{enumerate}
{\bf Proof}
\begin{enumerate}
\item As suggested in [3], to show the orbital stability of the
standing waves of (1.1), it suffices to prove that : $z_{c_1,c_2}
\neq \emptyset$ and any minimizing sequence\\
 $z_n = (z_{n,1},
z_{n,2}) \in \Sigma_{\mathbb{C}}(\mathbb{R}^2) \times
\Sigma_{\mathbb{C}}(\mathbb{R}^2)$ such that $\|z_{n,1}\|_2
\rightarrow c_1$   and  $|z_{n,2}\|_2 \rightarrow c_2$ and
$\widetilde{\hat{E}}(z_n)\rightarrow
\widetilde{\hat{I}_{c_1,c_2}}$\hfill(3.1)\\ is relatively compact in
$\sum_{\mathbb{C}}(\mathbb{R}^2) \times
\sum_{\mathbb{C}}(\mathbb{R}^2)$.\\
Let $z_n = (z_{n,1},z_{n,2})$ (with $z_{n,1} = (u_{n,1}, v_{n,1}),
z_{n,2} = (u_{n,2}, v_{n,2}) \subset \sum_{\mathbb{C}}(\mathbb{R}^2)
\times \sum_{\mathbb{C}}(\mathbb{R}^2)$ be a sequence such that
$\|z_{n,1}\|^2_2 \rightarrow c_1$ $\|z_{n,2}\|_2 \rightarrow c_2$
and $\widetilde{\hat{E}}(z_{n,1},z_{n,2}) \rightarrow
\widetilde{\hat{I}}_{c_1,c_2}$.
\end{enumerate}
Our first goal is to prove that $\{z_n\}$ has a subsequence which is
convergent in $\sum_{\mathbb{C}}(\mathbb{R}^2) \times
\sum_{\mathbb{C}}(\mathbb{R}^2)$.\\
By Lemma 2.2, it can be immediately deduced that $\{z_n\}$ is
bounded in $\sum_{\mathbb{C}}(\mathbb{R}^2) \times
\sum_{\mathbb{C}}(\mathbb{R}^2)$, therefore passing to a
subsequence, one can suppose that :
$$u_{n,i} \rightharpoonup u_i \mbox{ and } v_{n,i} \rightharpoonup
v_i \mbox{ in } \sum (\mathbb{R}^N), 1 \leq i\leq 2. \eqno{(3.2)}$$
Now set $\rho_{n,i} = |z_{n,i}| = (u^2_{n,i} + v^2_{n,i})^{1/2}$.\\
It follows that $\{\rho_{n,i}\} \subset \sum(\mathbb{R}^2)$ and that
for all $n \in \mathbb{N}$ and $1 \leq i,j \leq 2$ :
$$\partial_j\rho_{n,i} = \left\{ \begin{array}{ll}
\displaystyle{\frac{u_{n,i}(x)\partial_j
u_{n,i}(x)+v_{n,i}(x)\partial_j
v_{n,i}(x)}{(u^2_{n,i}+v^2_{n,i})^{1/2}}} &\mbox{ if }
u^2_{n,i}+v^2_{n,i} > 0\\
0 &\mbox{ elsewhere }
\end{array}\right.$$
Thus
$$\begin{array}{ll} \widetilde{\hat{E}}(z_n) -
\widetilde{E}(\rho_n) &= \displaystyle{\frac{1}{2}} \{\|\nabla
z_n\|^2_2 - \|\nabla
|z_n|\|^2_2\} \\
&= \displaystyle{\frac{1}{2}} \{\|\nabla z_{n,1}\|_2^2 + \|\nabla
z_{n,2}\|^2_2 - |\nabla
\rho_{n,1}|^2_2 - |\nabla \rho_{n,2}|^2_2\}\\
&= \displaystyle{\frac{1}{2}} \{|\nabla v_{n,1}|^2_2 + |\nabla
v_{n,1}|^2_2 + |\nabla u_{n,2}|^2_2 + |\nabla u_{n,2}|^2_2 - |\nabla
\rho_{n,1}|^2_2 - |\nabla \rho_{n,2}|^2_2\}\\
&= \displaystyle{\frac{1}{2} \sum^2_{i=1}
\int_{\{u^2_{n,i}+v^2_{n,i}> 0\}}} \sum^2_{j=1}
(\frac{u_{n,i}\partial_jv_{n,i}-v_{n,i}\partial v_{n,i})}
{u^2_{n,i}+v^2_{n,i}} \geq 0
\end{array}\eqno{(3.3)}$$
Hence $\widetilde{\hat{I}}_{c_1,c_2} =
\displaystyle{\lim_{n\rightarrow \infty}} \widetilde{\hat{E}}(z_n)
\geq \limsup \widetilde{E}(\rho_n)$.\\
Taking into account that
$$\|z_{n,i}\|^2_2 = |\rho_{n,i}|^2_2 = c^2_{n,i} \rightarrow c^2_i\quad
\forall\; 1 \leq i \leq 2,\eqno{(3.4)}$$
 Thus using Lemma 2.2 4), we
obtain that :
$$\liminf \widetilde{E}(\rho_n) \geq \liminf \widetilde{I}_{c_{n,1},c_{n,2}}
 \geq \widetilde{I}_{c_1,c_2}\geq
\widetilde{\hat{I}}_{c_1,c_2},$$ and hence
$$\lim_{n\rightarrow \infty} \widetilde{E}(\rho_n) = \lim_{n\rightarrow \infty}
\widetilde{\hat{E}}(z_n) = \widetilde{\hat{I}}_{c_1,c_2} =
\widetilde{I}_{c_1,c_2}.\eqno{(3.5)}$$ On the other hand (3.3)
implies that for any $1 \leq i \leq 2$, we have :
$$\lim_{n\rightarrow \infty} \int_{\mathbb{R}^2}
|\nabla u_{n,i}|^2 + |\nabla v_{n,i}|^2 -
|\nabla(u^2_{n,i}+v^2_{n,i})^{1/2}|^2 dx = 0.\eqno{(3.6)}$$
(3.2)
together with (3.6) imply that $\forall\; 1 \leq i \leq 2$.
$$\lim_{n\rightarrow \infty}
\int_{\mathbb{R}^2}|\nabla u_{n,i}|^2 + |\nabla v_{n,i}|^2 dx =
\lim_{n\rightarrow \infty} \int_{\mathbb{R}^2} |\nabla
(u^2_{n,i}+v^2_{n,i})^{1/2}|^2\eqno{(3.7)}$$ which  is equivalent to
say that
$$\lim_{n\rightarrow \infty}\|\nabla z_n\|^2_2 = \lim_{n\rightarrow \infty}
\|\nabla |z_n|\|^2_2.\eqno{(3.8)}$$
Now using (3.4), (3.5) and
Remark 1, $\rho_n = (\rho_{n,1}, \rho_{n,2})$ is relatively compact
in $\sum (\mathbb{R}^2) \times \sum (\mathbb{R}^2)$. Thus, there
exist $\rho_1, \rho_2 \in \sum (\mathbb{R}^2)$ such that :
$$\left\{ \begin{array}{l}
(u^2_{n,i}+v^2_{n,i})^{1/2} \mbox{ converges to } \rho_j \mbox{ in }
\sum(\mathbb{R}^2) : \;\forall\; 1 \leq i \leq 2\\
|\rho_j|_2 = c_j \mbox{ and }\\
\widetilde{E} (\rho_1,\rho_2) = \widetilde{I}_{c_1,c_2}
\end{array}\right.\eqno{(3.9)}$$
Let us first prove that $\rho_i = |z_i| = (u^2_i+v^2_i)^{1/2}$ ;
($u_i$
and $v_i$ are given in (3.2)).\\
By (3.2), we know that $u_{n,i} \rightarrow u_i$ and $v_{n,i}
\rightarrow v_i$ in $L^2(B(0,R))$, and we can easily see that :
$$[(u^2_{n,i} + v^2_{n,i})^{1/2} - (u^2_i + v^2_i]^2
\leq |u_{n,i}-u_i|^2+ |v_{i,n}-v_i|^2,$$
therefore
$$(u^2_{n,i}+v^2_{n,i})^{1/2} \rightarrow (u^2_i+v^2_i)^{1/2}
\mbox{ in } L^2(B(0,r))\quad \forall\; R > 0.$$ But
$(u^2_{n,i}+v^2_{n,i})^{1/2}\rightarrow \rho_i$ in $L^2(B(0,R)$,
this certainly implies that
$|z_i| = \rho_i\;\; \forall\; 1 \leq i \leq 2$.\\
On the other hand $\|z_{n,i}\|_2 = \|z_{n,i}\|_2 \rightarrow c_i =
\|z_i\| = ||z_i||_2$.\\
Therefore the proof of the first part of Theorem 3.1 is complete if
we show that $$\lim_{n\rightarrow \infty} \|\nabla z_{n,i}\|^2_2
\rightarrow \|\nabla z_i\|^2_2\;\quad \forall\; 1 \leq i \leq 2.$$
From (3.6), we have that
$$\lim_{n\rightarrow \infty}
\|\nabla z_{n,i}\|_2 = \lim_{n\rightarrow \infty} | \nabla
|z_{n,i}||_2 \quad \mbox{ and } \quad \lim_{n\rightarrow + \infty}
|\nabla|z_{n,i}||_2 = |\nabla|z_i\|^2_2.$$ Hence by the lower
semi-continuity of $|\;|_2$, we have :
$$\|\nabla z_i\|^2_2 \leq \lim |\nabla z_{n,i}\|^2_2 = \lim |\nabla |z_{n,i}\|^2_2 =
|\nabla |z_i||^2_2\eqno{(3.10)}$$ Finally, replacing $z_{n,i}$ by
$z_i$ in (3.3), we see that :
$$\|\nabla z_i\|^2_2 \geq |\nabla |z_i\|^2_2\quad \forall\; 1 \leq i \leq 2.$$
Now using (3.2), we know that $z_{n,i} \rightarrow z_i$ in
$\sum_{\mathbb{C}}(\mathbb{R}^N)$. Thus $z_{n,i} \rightarrow z_i$ in
$\sum_{\mathbb{C}}(\mathbb{R}^N)$ $\forall\; 1 \leq i \leq 2$.\\
{\bf Proof of 2) } Let $z = (z_1,z_2) \in z_{c_1,c_2}$ with $z_1 =
(u_1,v_1)$ and $z_2 = (u_2,v_2)$.

Let $\rho_1 = (u^2_1+v^2_2)^{1/2}$ and $\rho_2 =
(u^2_2+v_2^2)^{1/2}$. By the latter, we certainly have that
 $$\forall\; 1 \leq i \leq 2\;\;\forall\; 1\leq j \leq \;
 \int_{\mathbb{R}}\left(\frac{u_i\partial_jv_i-v_i\partial_ju_i}
 {u^2_i+v^2_i}\right)^2
 dx = 0\eqno{(3.11)}$$
 On the other hand $\widetilde{\hat{E}}(z_1,z_2) =
 \widetilde{\hat{I}}_{c_1,c_2}$, which implies that there exists a
 Lagrange multiplier $\alpha \in \mathbb{C}$ such that
 $$\widetilde{\hat{E}}(z)\xi = \frac{\alpha^n}{2} \sum^2_{i=1}
 z_i\bar{\xi}_i +  \overline{\xi z_i} \mbox{ for all } \xi \in
 \mathbb{C} \times \mathbb{C}.$$
 By elementary regularity theory and maximum principle, we can prove
 that $u_i$ and $v_i \in C^1(\mathbb{R}^2) \cap \sum (\mathbb{R}^2)$
 and $\rho > 0$.

 Set $\Omega = \{x \in \mathbb{R}^2 : u_i (x) = 0\}$ then
 $\Omega$ is closed since $u_i$ is continuous. Let us prove that
 it is also open.\\
 Let $x \in \Omega$. Using the fact that $v_i(x) > 0$, we can find a
 ball $B$ centered in $x_0$ such that $v_i(x) \neq 0$ for any $x \in
 B$.\\
 Thus for $x \in B$
 $$(\frac{v_i \partial_jv_i-v_i\partial_jv_i)^2}{u^2_i+v^2_i} =
 [\partial_j(\frac{u_i}{v_i}]^2
 \frac{v^4_i}{u^2_i+v^2_i} \mbox{ for } 1 \leq i,j \leq 2.$$
 This implies that
 $$\int_B|\nabla (\frac{u_i}{v_i})|^2\; \frac{v^4_i}{u^2_i+v^2_i} = 0 .\eqno{(3.8)}$$
 Hence $\nabla(\frac{u_i}{v_i}) = 0$ on $B \Rightarrow \;\exists\;
 C$ such that $\frac{u_i}{v_i} = C$ on $B$. Since $x_0 \in B \Rightarrow C \equiv
 0$.\\
 Therefore $\Omega$ is also an open set of $\mathbb{R}^N$. Hence we
 have proved that for $1 \leq i \leq 2$, these are two alternatives
 :
 \begin{enumerate}
 \item $u_i \equiv 0$ or $u_i > 0$ or $\mathbb{R}^2$
 \item $v_i \equiv 0$ or $v_i > 0$ or $\mathbb{R}^2$.
 \end{enumerate}
 Now let $z_i = e^{i\sigma_i}w_i$, $\sigma_i \in \mathbb{R}$, $w_i \in
 W_{c_1,c_2}$.
 Thus $|z_i|_2 = c_i$ and $\widetilde{E}(z_1,z_2) =
 \hat{I}_{c_1,c_2}$.
 Then $\{e^{i\sigma_i}w_i : \sigma_i \in \mathbb{R},
 w_i \in W_{c_1,c_2}\}\subset Z_{c_1,c_2}$.\\
 Conversely for $z_i = (u_i,v_i)$ such that $(z_1,z_2) \in
 Z_{c_1,c_2}$, set $w_i = |z_i|$. Then $\widetilde{\hat{E}}(z_1,z_2)
  = \widetilde{E}(w_1,w_2) =
 \widetilde{\hat{I}}_{c_1,c_2} = \hat{I}_{c_1,c_2}$ and $(w_1,w_2) \in
 W_{c_1,c_2}$.\\
 We now have four possible alternatives. We will discuss one in
 details, the three others can be shown following exactly the same
 ideas.

 Suppose that $v_1$ and $v_2 \neq 0$ for all $x \in \mathbb{R}^2$.

 In this case, it follows that $\nabla (\frac{u_i}{v_i}) = 0$ on
 $\mathbb{R}^2$.\\
 Thus we can find 2 constants $K_1, K_2 \in \mathbb{R}$ such that
 $$u_1 \equiv K_1v_1\quad \mbox{ and }\quad u_2 \equiv K_2 v_2.$$
 Therefore $w_1 = (K_1+i)v_1$ and $w_1 = |K_1+i||v_1|$.\\
 Let $\theta_1 \in \mathbb{R}$ such that $K_1 +i = |K_1+i|e^{\theta_1}$
 and let $\varphi_1 = 0$ if $v_i > 0$ and $\varphi_1 = \pi$ if $v_1 <
 0$ on $\mathbb{R}^2$. Setting $\sigma_1 = \theta_1 + \varphi_1$, $z_1 =
 (K_1+i)v_1 = |K_1+i| e^{\theta_1}|v_1|e^{i\varphi_1} =
 w_1 e^{i\sigma_1}$.\\
 Similarly $z_2 = w_2 e^{i\sigma_2}$ with $w = (w_1,w_2)$.
 \section*{References}
 \begin{enumerate}
 \item {\bf P. Antonelli, Rada Maria Weishaeupl} : Asymptotic
 behaviour of nonlinear Sch\"odinger systems with Linear coupling,
 preprint.
 \item {\bf W. Bao, Y. Cai} : Ground states of two-component
 Bose-Einstein condensates with an internal Josephson function. East
 Asian Journal of applied Mathematics, Vol 1, No 1 49-81, 2011.
 \item {T. Cazenave, P.L.Lions} Orbital stability of standing waves
 for some Schr\"odinger equations. Comm math Phys, 85, p. 549-561
 (982).
 \item {\bf H. Hajaiej} : Cases of equality and strict inequality in
 the extended Hardy-Littlewood inequalities. Proc. Roy. Pro
 Eduburgh, 135 A, (2005), 643-661.
 \item {\bf H. Hajaiej} : Extended Hardy-Littlewood iqualities and
 some applications, Trans. Amer. Math. Soc, Vol 357, No 12, 20054 pp
 4885-4896.
 \item {\bf H. Hajaiej, C.A. Stuart} : On the variational approach
 to the stability of standing waves for nonlinear Schr\"odinger
 equations adv. Nonlinear Studies, 4 (2004), 469-501.
 \item {\bf A. Jungel, Rada. Maria Weishaeupl} Blow-up in two
 component nonlinear Schr\"odinger systems with  an external driven
 field. M3AS (in press).
 \item {\bf O. Kavian, F.B.Weissler} : The pseudo-conformally
 invariant nonlinear Schr\"odinger equation, Michigan Math J 41,
 151-173.
 \end{enumerate}
\end{document}